\documentclass[11pt]{article}
\usepackage{amsmath, mathtools, amssymb, amsthm, amsfonts}
\usepackage[numbers,comma,sort&compress]{natbib}
\usepackage{graphicx}
\usepackage{color}
\usepackage{caption}
\usepackage{subcaption}
\usepackage{float}
\usepackage{url}
\usepackage{authblk}

\date{}
\title{An Upper Bound Estimate and Stability for the Global Error of Numerical Integration Using Double Exponential Transformation  }
\author[1]{Arezoo  Khatibi}
\author[2]{Omid Khatibi\footnote{Omid Khatibi, correspond email: \protect\url{md_khn@yahoo.com,a1449933@unet.univie.ac.at}, Faculty of Mathematics,  University of Vienna,  Oskar Morgenstern platz 1 1090, Vienna Austria.}}
\affil[1]{University of Kashan,Kashan,Iran}
\affil[2]{University of Vienna,Vienna,Austria}

\begin{document}
\maketitle

\begin{abstract}
The double exponential formula  was introduced  for calculating  definite  integrals   with  singular point oscillation functions and Fourier integral. The double exponential  transformation is not only useful  for numerical  computations  but  it  is  also used  in  different  methods  of  Sinc  theory. In  this paper  we  give an upper bound estimate for the error of double exponential transformation. By  improving  integral  estimates  having  singular  final points, in theorem 1 we prove that the method  is  convergent and the rate of convergence is $\mathcal{O}(h^2)$  where h is a step size. Our main tool in the proof is DE formula in Sinc theory. The  advantage  of  our method  is  that  the  time  and  space  complexity  is  drastically  reduced. Furthermore, we discovered upper bound error in DE formula independent of N truncated number, as a matter of fact we proved stability. Numerical tests are presented to verify the theoretical predictions and confirm the convergence of the numerical solution.
 \end{abstract}
\section{Introduction}
The double exponential transformations (DE) is used for evaluation of integrals of an analytic function has end point singularity. Our innovation has discovered the upper bound error in DE formula independent of N truncated number, therefore we proved stability. First we introduce double exponential transformation method. The main idea in this method change the variable by a function $\phi$ which transforms the interval of integration change to $(-\infty,+\infty)$. In section 1 we had provided background introduced the results of reference [1] and in section 3 we presents our main results in theorem 1 and provide numerical test to verify our result. Tanaka and Sugihara make full study of function classes for successful DE-Sinc approximations [2]. In [3] Ooura expanded a double exponential formula for the Fourier transform. Ooura studied animt-type quadrature formula with the same asymptotic performance as the DE-formula [4]. Several scientist studied DE-formula [5-10]. Stenger used DE Formula in Sinc approximations [9].
\section{Double Exponential Transformation}

Consider the following integral 
\begin{equation}\label{eq:int}
\int_a^b\!f(x)\,\mathrm{d}x,
\end{equation}
where the interval $(a,b)$ is infinite or half infinite and the function under integral is analytic  on $(a,b)$ and perhaps has singular  point  in $x=a$   or $x=b$ or both, now consider the  change  of  variables below [11]:
\begin{equation}\label{eq:varchange}
x=\phi(t),a=\phi(-\infty),b=\phi(+\infty),
\end{equation} 
 where $\phi$ is analytic on $(-\infty,+\infty)$ and 
 \begin{equation}\label{eq:int}
I=\int_a^b\!f(\phi(t))\phi\prime(t)\,\mathrm{d}t.
\end{equation}  
Such that after the change of variable the integrand decays double exponential:
\begin{equation}\label{eq:int}
\mid f(\phi(t))\phi\prime(t)\mid\approx e^{-c e^ {\mid t \mid}}, \mid t \mid\rightarrow\infty, c> 0. 
\end{equation}
By using   the trapezoidal formula with mesh size h on (3) we have 
\begin{equation}\label{eq:trapez}
I_{h}=h\sum_{k=-\infty}^{+\infty}f(\phi(kh))\phi\prime(kh). 
\end{equation}
The above infinite summation is truncated from $k=〖-N〗^-$ to $k =N^+$ in computing \eqref{eq:trapez}  therefore we conclude,
\begin{equation}\label{eq:int}
I_{h}^{(N)}=h\sum_{k=-N^-}^{N^+}f(\phi(kh))\phi\prime(kh),N=N^+ +N^- +1 
\end{equation}

Here the N+ 1 relates to zero point, for example, $N=50$ we have $N^-=24, N^+ =25 $ finally 50=25+24+1.
What gives us the authority to truncate n numbers   from the   infinite summation? 
 Indeed $e^{(-ce^{\mid t \mid} )}$ quickly   approaches zero.
 An  analysis of  the   graph $e^{(-ce^{\mid t \mid} )}$  for small values  of  t  and  (4)  with c=2 is  as  follows  Graph A shows the double exponential function.
 
    \begin{figure}[htbp]
 \centering
    \begin{subfigure}[b]{0.3\textwidth}
        \includegraphics[width=\textwidth]{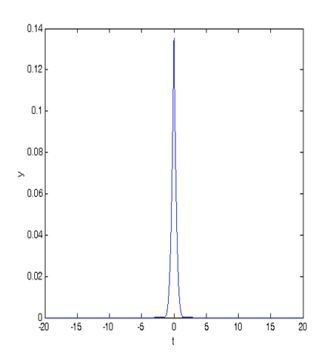}
        \caption{Graph A}
        \label{Graph A}
    \end{subfigure}
    \begin{subfigure}[b]{0.31\textwidth}
        \includegraphics[width=\textwidth]{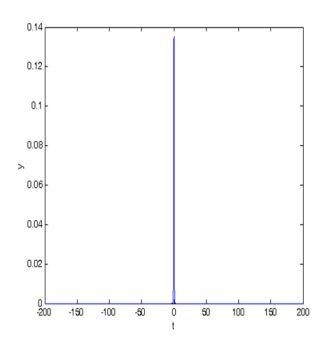}
        \caption{Graph B}
        \label{Graph B}
    \end{subfigure}
\end{figure}

Graph B  clearly shows  the expansions   from both side close to axis x  and  error  growth  is slow and arrives  at zero  as   we use  the double  exponential  formula   and  abbreviate  it as DE.  
For the integral over (-1, 1).
\begin{equation}\label{eq:int}
\int_ {-1}^ 1\!f(x)\,\mathrm{d}x,
\end{equation}
let
\begin{equation}\label{eq:int}
\phi(t) = \tanh (\frac{\pi}{2 }\sinh(t)),
\end{equation}
the transformation $ x= \phi(t) = \tanh (\pi/2 \sinh(t)) $ satisfied and substitute in (5) we have
\begin{equation}\label{eq:int}
\phi(t)\prime= \tanh \prime (\pi/2\sinh(t))\pi/2\ cosh(t)= \left(\frac {1}{\cosh^2(\pi/2\sinh(t))}\right)\pi/2\cosh(t).
\end{equation}
After using double exponential formula we will obtain
\begin{equation}\label{eq:int}
I_{h}^{(N)}=h\sum_{k=-N^-}^{N^+}f(\tanh(\pi/2\sinh(kh))\frac{ \pi/2 \cosh(kh)}{\cosh^2 (\frac{\pi \sinh(kh)}{2})}. \text{[11]} 
\end{equation}
In the next section, we compute the error in (5)  with trapezoidal rule and mesh size h.
\section{Upper Bound Error, Convergence Rate and Stability of D.E Function}\newtheorem{theorem}{Theorem}
\begin{theorem}

$\forall h>0, h$ is step size, $c>0$ and $c$ is constant, let \begin{equation}
I=\int_{-\infty}^\infty\!f(\phi(t))\phi\prime(t)\,\mathrm{d}t,
\end{equation}
class function of f is valid for analytic function has a singularity   in first point or end point or both,
$\phi$ provided by 
\begin{equation}
\phi(t) = \tanh (\frac{\pi}{2 }\sinh(t)),
\end{equation} 
such that
\begin{equation}
 |F(t)|=|f(\phi(t))\phi\prime(t))|\approx  e^{(-ce^{|t|})},
 \end{equation}
let 
\begin{equation}
I_{h}^{(N)}=h\sum_{k=-N^-}^{N^+}f(\tanh(\pi/2\sinh(kh))\frac{ \pi/2 \cosh(kh)}{\cosh^2 (\frac{\pi \sinh(kh)}{2})},
 \end{equation}
 then the error estimate 
 \begin{equation}
 |I-I_{h}^{(N)}|  \leq h^2(\frac{e^{-4}(1+c)}{3}e^{\frac{-c}{2}}\frac{1}{12}(c+c^2)),
 \end{equation}
  holds true.

\begin{proof}

 The absolute error in  $k^{th}$ open interval [kh,(k+1)h] using trapezoidal rule is,
\begin{equation}\label{eq:int}
E_{k}=\frac{1}{12} h^3 M_{k}, 
\end{equation}
 \begin{equation}\label{eq:int}
M_{k}=\sup\mid F''(X) \mid,        
 kh\leq x \leq(k+1)h. 
\end{equation}
In the next relations, the function decreases, therefore $M_{k}=F''(kh)$ ,
\begin{equation}\label{eq:int}
F (t)=f(\phi(t))\phi\prime(t)), 
\end{equation}
thus,
\begin{equation}\label{eq:int}
|F(t)|\approx  e^{(-ce^{|t|})} 
, 
\end{equation}
\begin{equation}\label{eq:int}
\acute{F}(t)\approx-ce^{|t|} e^{(-ce^{|t|})} 
, 
\end{equation}
such that,
\begin{equation}\label{eq:int}
|\acute{F}(t)|\approx ce^{|t|} e^{(-ce^{|t|})} 
, 
\end{equation}
we obtain
\begin{equation}\label{eq:int}
|{F^{''}}(t)|\approx ce^{|t|} e^{(-ce^{|t|})}+c^2 e^{2|t|}
e^{-ce^{|t|})}.
\end{equation}
On the other hand, we know  that $e^{|t|}\leq e^{2|t|}$   and by 
 using the  triangle inequality
$|a+b|\leq |a|+|b|$  we  have:
\begin{equation}\label{eq:int}
|{F^{''}}(t)| \leq ce^{2|t|} e^{(-ce^{|t|})}+c^2 e^{2|t|}e^{-ce^{|t|}},
\end{equation}
\begin{equation}\label{eq:cpluscsquar}
\leq (c+c^2)e^{2|t|} e^{-Ce^{|t|}}.
\end{equation}
We know that  if $x\rightarrow\infty$
 then $e^{x}\rightarrow+\infty$ 
and if $x\rightarrow-\infty$ then  $e^{x}\rightarrow0$; if $x\rightarrow\pm\infty$
 then $|x|\rightarrow\infty$ if  $c>0$  then $e^{-c|x|}\rightarrow 0$.  
Note that  in \eqref{eq:cpluscsquar}, $2t-ce^t$  decreases  because 
$\frac{d}{dt}(2t-ce^t)=2-ce^t$ and 
for $t>ln\frac{2}{c}$ the derivation is negative and $kh\leq x \leq(k+1)h,$ therefore,
 $\sup\mid F''(X) \mid\leq(c+c^2)e^{2kh-ce^{kh}},$ the function is also decreasing. 
By  using  \eqref{eq:cpluscsquar}  we have \eqref{eq:sigmaek}.
We truncate  the  infinite  summation  from  $-N^-$ to $N^+$:
\begin{equation}\label{eq:sigmaek}
|I-I_{h}^{(N)}|\leq\sum_{k=-N^-}^{N^+}E_{k}\leq \frac{1}{12}
h^3\sum_{k=-N^-}^{N^+}M_{k}\leq \frac{1}{12}
h^3 (c+c^2)\sum_{k=-N^-}^{N^+}e^{2kh-ce^{kh}}.
\end{equation}

\paragraph{ Lemma (3.1):}There exists $k_0$ such that if $k>k_0$ 
then,\begin{equation}\label{eq:eminesfractionctwo}
e^{2kh-ce^{kh}}<e^{-\frac{c}{2} e^{kh}}.
\end{equation}
\begin{proof}

 Multiplying $e^{(ce^{kh)}}$in both sides of \eqref{eq:eminesfractionctwo} we obtain
\begin{equation}
e^{2kh}<e^{\frac{c}{2}e^{kh}}.                                     
\end{equation}
This  is enough  to show that, \begin{equation}\label{eq:inequalityorbetween} 2kh<\frac{c}{2}e^{kh}\quad \text{or}\quad  
\frac{4}{c}kh<e^{kh},\qquad \qquad\qquad\qquad \quad \qquad\qquad \quad
\end{equation}

so, we  need the  following  lemma  with t= kh and  which means a =4/c  can be used in next  lemma .
\end{proof}
\paragraph{ Lemma (3.2):}
 For $a>0$ and for $t> 2a $ such 
 we have $e^t>at$.
    
\begin{proof}
\begin{equation} 
e^t=1+\frac{t}{1!}+\frac{t^2}{2!}+\frac{t^3}{3!}+\cdots>1+t+\frac{t^2}{2},  
    t>0, 
\end{equation}
for this we prove that\begin{equation}1+t+\frac{t^2}{2}>at,\end{equation} 
\begin{equation} 
1+(1-a)t+\frac{t^2}{2}>0, 
\end{equation}
this yields
\begin{equation} 
t>\frac{a-1+\sqrt{(a-1)^2-2}}{1}, 
\end{equation}
also
\begin{equation} 
t>\frac{a-1-\sqrt{(a-1)^2-2}}{1}, 
\end{equation}
\text{ the discriminant of the quadratic equation } 
\begin{equation}\label{eq:mobinmoadele}
1+(1-a)t+\frac{t^2}{2}=0, 
\end{equation} 
\text{is }
\begin{equation}\label{eq:sefrmobinmoadele}
\Delta=(1-a)^2-2=a^2-2a-1=0,  
\end{equation}
if $\Delta\geq 0$ we  put  $t_.$  equal  to the max root in quadratic 
equation \eqref{eq:mobinmoadele},
\begin{equation} 
t_.=a-1+\sqrt{(a-1)^2-2}. 
\end{equation}
Then  the statement  is satisfied. And  note  that, 
\begin{equation}
 t_.< a+\sqrt{a^2}=2a.
\end{equation}
If $\Delta<0$  that is  the roots are  complex  and the coefficient of $t^2$ 
is positive then  for all of  t  inequality is true $.$ Hence  we  can  set 
$t_.=0$   
and the proof is complete.
If $\Delta\geq 0$  we obtain:
\begin{equation}
 t>2 a, a =\frac{4}{c}, t=kh \Rightarrow kh>\frac{8}{c}. 
\end{equation}
\begin{equation}
 t<2 a, a =\frac{4}{c}, t=kh \Rightarrow kh<\frac{8}{c}. 
\end{equation}
\end{proof}
In \eqref{eq:sigmaek}  we   get:
\begin{equation}\label{eq:globalerrorghablcase1}\text{Global  Error}\leq \frac{1}{12} h^3 (c+c^2) 
\left( \sum_{k=[\frac{8}{ch}]}^{N^+}e^{2kh-ce^{kh}}+ 
\sum_{k=-N^-}^{[\frac{8}{ch}]}e^{2kh-ce^{kh}}\right)
\end{equation} 
Now  we  consider  both  sums  separately. \\
\textbf{Case 1: }
\begin{equation}\frac{1}{12}h^3 (c+c^2) \sum_{kh>\frac{8}{c}}e^{2kh-ce^{kh}}.
\end{equation}
\begin{proof}

Put  \eqref{eq:inequalityorbetween}  in   case1   and  use  lemma (3.2) to set,  
\begin{equation}
\leq \frac{1}{12} h^3 (c+c^2)
\sum_{k=[\frac{8}{ch}]}^{N^+}e^{-\frac{c}{2} e^{kh},}
\end{equation}
On the other hand we know that  for 
$t>0$ we have, \begin{equation}
e^t>1+t,
\end{equation}     
 therefore,\begin{equation}
e^{|kh|} >1+|kh|,
\end{equation}
\begin{equation}
\sum_{[\frac{8}{ch}]}^{N^+}e^{-\frac{c}{2} e^{kh}} 
\leq\sum_{k=[\frac{8}{ch}]}^{N^+}e^{-\frac{c}{2}-
\frac{c|kh|}{2}}
\leq e^{-\frac{c}{2}}
\sum_{k=[\frac{8}{ch}]}^{N^+}e^{-\frac{c}{2}|kh|} 
\leq e^{-\frac{c}{2}} 
\sum_{k=[\frac{8}{ch}]}^{N^+}e^{-\frac{ch}{2}|k|},
\end{equation}
by using the geometric series we obtain:
\begin{equation}
\sum_{k=k_.}^n\frac{a^{k_.}-a^{n+1}}{1-a}<\frac{a^{k_.}}{1-a},
\end{equation}
such that
\begin{equation}
a= e^{-\frac{ch}{2}}, k_.=[\frac{8}{ch}],
\end{equation}
\begin{equation}
\leq e^{-\frac{c}{2}}
\frac{(e^{-\frac{ ch}{2}})^{\frac{8}{ch}}}
{1-e^{-\frac{ch}{2}}},
\end{equation}
moreover,
\begin{equation}
\leq e^{\frac{-c}{2}}\frac{e^{-4}}{1-e^{-\frac{ch}{2}}},
\end{equation}
therefore   \eqref{eq:sefrmobinmoadele}  gives the following inequality:
\begin{equation}
\text{Error}\leq \frac{c+c^2}{12}e^{\frac{-c}{2}}
\frac{h^3 e^{-4}}{1-e^{\frac{-ch}{2}}}.
\end{equation}

 \begin{equation}\label{eq:soratlemsese}  \textbf{Lemma (3.3)}:\text{There  exists $h_0$  such that if } 
h< h_0:e^{\frac{-ch}{2}}\leq 1-\frac{ch}{4}.
\end{equation}
\begin{proof}

\begin{equation}
f(h)= e^{\frac{-ch}{2}}-(1-\frac{ch}{4}),\quad\quad\quad\quad\quad\quad\quad\quad\quad\quad\quad\quad\quad\quad
\end{equation}
we have $f(0)=1-1=0$, and
\begin{equation}
f'(h)=\frac{-c}{2} e^{\frac{-ch}{2}}+\frac{c}{4},
\end{equation}
because  $\lim_{h\rightarrow 0 }e^{\frac{-ch}{2}}=1$
then
\begin{equation}
\lim_{h\rightarrow 0 }f'(h)=\frac{-c}{2} +\frac{c}{4}=\frac{-c}{4},
\end{equation}
there  exists  $h_0>0$ so that if $h<h_0$ then $f'(h)<0$  
and $f$
 decreases  because $f(0)=0$  and  for $h<h_0$,
 $f(h)\leq 0$;  hence $e^{\frac{-ch}{2}}\leq 1-ch/4$.

\end{proof}
\begin{equation}
1-e^{\frac{-ch}{2}}\geq \frac{ch}{4},
\end{equation}
we infer that 
\begin{equation}
\frac{1}{1-e^{\frac{-ch}{2}}}\leq \frac{1}{\frac{ch}{4}}\leq\frac{4}{ch},
\end{equation}
and
\begin{equation}
\leq \frac{c+c^2}{12}e^{\frac{-c}{2}}
\frac{h^3}{1-(1-\frac{ch}{4})},
\end{equation}
thus,
\begin{equation}
\leq \frac{c+c^2}{12} e^{\frac{-c}{2}}\frac{h^3\times 4}{ch}(e^{-4}),
\end{equation}
clearly
\begin{equation}
\leq \frac{(1+c)}{3}
 e^{\frac{-c}{2}} (e^{-4}) h^2,
\end{equation}
indeed
\begin{equation}
\leq \frac{e^{-4}(1+c)}{3} e^{\frac{-c}{2}}h^2,
\end{equation}
hence the case 1 is completed.
\end{proof}
\textbf{Case 2:}
\begin{equation}
\frac{1}{12} h^3 (c+c^2)\sum_{kh<\frac{8}{c}}e^{2kh-ce^{kh}}.
\end{equation}
\begin{proof}
\begin{equation}
2kh-ce^{kh} < 2kh,
\end{equation}
\begin{equation}
\sum_{kh<\frac{8}{c}}e^{2kh-ce^{kh}}<\sum_
{k=-N^{-}}^{[\frac{8}{ch}]}e^{2kh}<
\sum_{k=[\frac{8}{ch}]}^{-\infty}(e^{2h})^k=
\sum_{k=-[\frac{8}{ch}]}^{\infty}(e^{-2h})^k,
\end{equation}
we know
\begin{equation}
\sum_{k=k_.}^\infty a^k =\frac{a^{k_.}}{1-a},
\end{equation}
then
\begin{equation}
\sum_{k=-[\frac{8}{ch}]}^\infty(e^{-2h})^k=
\frac{(e^{-2h})^{-[\frac{8}{ch}]}}{1-e^{-2h}},
\end{equation}
so that
\begin{equation}
[\frac{8}{ch}]\leq \frac{8}{ch}, -[\frac{8}{ch}]\geq -\frac{8}{ch},
e^{-2h}<1,
\end{equation}
therefore,
\begin{equation}
\leq \frac{(e^{-2h})^{\frac{-8}{ch}}}{1-e^{-2h}},
\end{equation}
\end{proof}
in   \eqref{eq:globalerrorghablcase1}  we obtain:
\begin{equation}
|I-I_{h}^{(N)}| \leq \frac{e^{-4}(1+c)}{3} e^{\frac{-c}{2}} h^2+\frac{1}{12} h^3 
(c+c^2)\frac{(e^{-2h})^{\frac{-8}{ch}}}{1-e^{-2h}},
\end{equation}
in  \eqref{eq:soratlemsese}  consider,
\begin{equation}
\frac{ch}{4}=h
\Rightarrow e^{-2h}\leq 1-h\text{.}
\end{equation}
This yildes
\begin{equation}|I-I_{h}^{(N)}|   
\leq \frac{e^{-4}(1+c)}{3} e^{\frac{-c}{2}}h^2+\frac{1}{12} h^3 
(c+c^2)\frac{e^{\frac{16}{c}}}{1-(1-h)}\text{,}
\end{equation}
\begin{equation}|I-I_{h}^{(N)}|   
\leq(h^2\frac{e^{-4}(1+c)}{3}e^{\frac{-c}{2}} +\frac{1}{12}  
(c+c^2))\text{.}
\end{equation}
Now  let $-N^-,N^+\rightarrow +\infty$. So we have now proved  
theorem 1.\end{proof}\end{theorem}
\textbf{Remark:} Note that the upper bound error of theorem 1 is independent of N truncated number that means we proved \textbf{stability}.
If $h\rightarrow 0$ error approaches zero  and integrand method is {\bf convergence}
  and 
accuracy is $\mathcal{O}(h^2)$. \\
By using this method we improve the space complexity and 
time complexity and we also increase the rate of convergence .
In fact, we have found upper bound for global error {\bf (u.b.g.e)}.
\section{ Numerical Result}
In the  next example, utilising  numerical computation and Maple  software to show  the result. Assume
 \begin{equation}
 \quad I1=\int_0^1\exp(20(x-1))\sin(256x)dx\text{.}
\end{equation}
The  absolute  error  tolerance  is $10^{-8}$ where  $N$  is the 
number of function  evaluations and abs. error  is the actual absolute  error 
of  the result and  it gives an approximate value  which is correct up to 16 
significant digits.

\noindent
clc\\
clear all\\
tic\\
$h=1/129$; $c=2$;\\
GError=$h^2/3*(1+c)*(exp(-4-c/2)+1/4*c)$\\
toc\\
Results:

GError = $ 3.0451e-005$

Elapsed time is $0.000106$ seconds.
\paragraph{TABLE  1 comparison}
Comparison  of  the  efficiency  of {\bf DEFINT}
  and {\bf ubge}. The  absolute  error tolerance  is $10^{-8}$ where  $N$  is the number of function  
evaluations and abs. error  is the actual absolute  error of  the result and it gives an approximate value  which is correct up to 16 significant digits .

\begin{center}
\begin{tabular}{cccc}\hline
\qquad\qquad\qquad\qquad DEFINT\\\hline
INTEGRAL & N & abs. error[1] & ubge\\\hline
$I_1$ &259   & $4.8\times 10^{-12}$ & $3\times 10^{-5}$\\\hline
\end{tabular}
\end{center}
Considering the performance  of convergence  rate of quadrature rules  by  
using D.E function, the present formula  is suitable for automatic integrator, sinc method, Fourier integral and increases the rate of convergence. It  
Improves the  space complexity and  the time complexity  and has an important 
role  in numerical  analysis.


\begin{thebibliography}{9}
\bibitem{M. Mori , M. Sugihara2001}
  M. Mori , M. Sugihara,
  \emph{\:The Double Exponential Transformation in Numerical  Analysis },
  Journal Of  Computational And Applied  Mathematics Vol.127,
  (2001), 287-299.
  \bibitem{K. Tanaka, M. Sugihara, K. Murata2009}
 K. Tanaka, M. Sugihara, K. Murata,
  \emph{\:Function Classes for Successful DE-Sinc Approximation  },
  Mathematics of Computation,Vol.78 ,
  (2009), 1553-1571.
  \bibitem{T. Ooura2005}
 T. Ooura,
  \emph{\:A Double Exponential Formula for Fourier Transform  },
   RIMS, Kyoto Univ., Vol. 41 ,
  (2005), 971-978.
  \bibitem{T. Ooura2008}
 T. Ooura,
  \emph{\:Animt-Type Quadrature Formula with the Same Asymptotic Performance as the De-Formula  },
   Journal of Computational and Applied Mathematics, Vol. 213 ,
  (2008), 232-239.
   \bibitem{K. Horiouchi, M.Sugihara1999}
 K. Horiouchi, M.Sugihara,
  \emph{\:Sinc - Galerkin  Method  with  The Double Exponential Transformation for the Two Point Boundary Problems  },
   Technical Report, Dept. of Mathematical Engineering, University of Tokyo ,
  (1998), 99-105.
  \bibitem{M.Mori, A.Nurmuhammad, M.Muhammad2009}
 M.Mori, A.Nurmuhammad, M.Muhammad,
  \emph{\:De-Sinc Method for Second Order Singularly Perturbed Boundary Value Problems  },
   japan j.Indust.Appl.Math, Vol. 26 ,
  (2009), 41-63.
  \bibitem{M. Muhammad, A. Nurmuhammad, M. Mori2005}
 M. Muhammad, A. Nurmuhammad, M. Mori,
  \emph{\:Numerical  Solution of Integral Equations by Means of the Sinc collocation Method Based on the Double Exponential Transformation   },
   Journal of Computational  and Applied Mathematics , Vol. 177 ,
  (2005), 269-286.
  \bibitem{M. Sugihara1997}
 M. Sugihara,
  \emph{\:Optimality of The Double Exponential Formula –Functional Analysis Approach   },
   Numer Math, Vol. 75 ,
  (1997).
  \bibitem{F. Stenger1993}
 F. Stenger,
  \emph{\:Numerical Methods Based on Sinc and Analytic Functions },
     
  (1993).
  \bibitem{T. Koshihara, M. Sugihara1996}
T. Koshihara, M. Sugihara,
  \emph{\:Optimality of The Double Exponential Formula –Functional Analysis Approach   },
   Numer Math, Vol. 75Proceedings Of The 1996 Annual Meeting Of The Japan Society For Industrial And Applied Mathematics ,
  (1996).
  
  \bibitem{O. Khatibi, A. Khatibi2017}
O. Khatibi, A. Khatibi,
 \emph{\:Criteria for the Application of Double Exponential
Transformation  },
arXive, Preprints.
 \end{thebibliography}
\end{document}